\documentclass[11pt]{article}  
\usepackage{amsmath,amsfonts,amssymb,color}
\usepackage{graphicx}
\usepackage[mathscr]{euscript}
 \usepackage{setspace}
\newcommand{\bs}{\boldsymbol}
\newcommand{\tx}{\textstyle}

\newcommand{\eps}{\epsilon}

\newcommand{\pr}{\prime}
\newcommand{\al}{\alpha}

\newcommand{\rw}{\rightarrow}

\newcommand{\R}{\mathbb R}

\begin{document}   
\title{Nonsmooth Invariant Manifolds in a Conceptual Climate Model}
\author{James Walsh\\Department of Mathematics, Oberlin College}
\date{May 8, 2017}

\maketitle

\vspace{.25in}
\noindent{\bf Abstract.} \ There is widespread agreement that ice sheets   flowed into the ocean in tropical latitudes at sea level during the Earth's past. Whether these extreme ice ages were snowball Earth events, with the entire surface covered in ice, or whether   ocean water remained ice free in regions about the equator, continues to be controversial. For the latter situation to occur, the effect of positive ice albedo feedback would have to be damped  to stabilize an advancing ice sheet shy of the equator. In this paper we analyze a conceptual model comprised of a zonally averaged surface temperature equation coupled to a dynamic ice line equation.  This difference equation model is aligned with the cold world of these great glacial episodes through an appropriately chosen albedo function. Using the spectral method, the analysis leads to a nonsmooth singular perturbation problem. The Hadamard graph transform method  is applied to prove the persistence of an invariant manifold, thereby providing insight into model behavior. A stable climate state with the ice line resting in tropical latitudes, but with open water about the equator, is shown to exist. Also presented are local smooth and nonsmooth bifurcations as parameters related to atmospheric CO$_2$ concentrations and the efficiency of meridional heat transport, respectively, are varied.

\bigskip



\section{Introduction} 
The theory of nonsmooth dynamical systems is playing an ever expanding role in the analysis of conceptual climate models. Nonsmooth models arising in the study of the overturning ocean circulation include \cite{julie},  \cite{stom}, \cite{verd} and \cite{wel}. Several works have appeared concerning nonsmooth models of sea ice loss in the Arctic  (\cite{eis}, \cite{eiswett}, \cite{hill}, for example).  Moreover, nonsmooth models  often appear when modeling the glacial cycles, where the lack of smoothness may be associated with the dynamic between ice sheets and atmospheric carbon dioxide \cite{anna} or the release of carbon dioxide from the deep ocean \cite{hogg}; with changes in Antarctic ice volume \cite{par}; or with changes in glacial ablation rates \cite{us}.

The present work focuses on surface albedo, which  plays a significant role in the evolution of a planet's climate. The term {\em albedo} refers to the extent to which the surface reflects incoming solar radiation (or {\em insolation}). Snow and ice, for example, have higher albedos than does ice free surface. 
Notably, ice albedo interacts with climate in a positive feedback loop. Were an existing ice sheet to expand,  the surface albedo would increase, thereby lowering temperature and triggering further ice sheet growth. On the other hand,  were an ice sheet to retreat, surface albedo would decrease, leading to greater absorption of insolation and an increase in temperature. This would lead to further reduction in the size of  the ice cover. (See \cite{noaa} for the effect the recent and ongoing loss of ice in the Arctic has had on Arctic temperatures, for example.)

There is widespread agreement that during two ice ages in the Neoproterozoic Era, roughly 715 and 630 million years (mya)  ago, respectively, ice sheets flowed into the ocean at sea level in tropical latitudes \cite{bender}. Many believe these extensive  glaciations were {\em snowball Earth} episodes, with ice covering the entire planet, triggered in part by positive ice albedo feedback (see \cite{hoff} and references therein).  There are others, however, who argue the advance of the ice sheet likely stopped shy of the equator, leaving a strip of open ocean water (see \cite{abb} and references therein). The actual extent of the ice cover during these  ice ages remains controversial \cite{bender}.

Conceptual climate models, focusing on notions of energy balance,  were introduced by M. Budyko \cite{bud} and W. Sellers \cite{sell} in 1969 to investigate positive ice albedo feedback. Having latitude as the sole spatial variable, each model focused on the effect ice albedo has on zonally averaged surface temperature.
 Sellers' model indicated  a reduction of the solar constant by 2-5\% was sufficient to trigger an ice age. Budyko's work suggested  that once the ice sheet extended down below a critical latitude, a snowball Earth event inexorably followed. The solar constant 600-700 mya was 6\% less than today's value, so that insolation absorbed at the surface was significantly lower during these times  relative to our present climate. The modeling of Sellers and Budyko, together with physical evidence \cite{hoff}, support the proposition  that the Neoproterozoic ice ages mentioned above were snowball Earth events.

There is also evidence, however, supporting the existence of open water near the equator during these ice ages \cite{abb}. For example, photosynthetic organisms and other life forms are know to have survived these glaciations.   A mechanism by which positive ice albedo feedback might be  sufficiently damped so as to ensure the advance of a large ice sheet stopped prior to reaching the equator was introduced in \cite{abb}. Abbot et al argued   that a stable climate state, with glaciers extending to the tropics but with open water encircling the equator, could exist given certain assumptions concerning the albedo, to be described below. With seasonal changes causing the open ocean to snake around the equator, this stable climate state was called a {\em Jormungand} state in \cite{abb}, after the sea serpent from Norse mythology.

We incorporate the  assumptions from \cite{abb} into an energy balance surface temperature model coupled with a dynamic ice line, assuming diffusive meridional  heat transport. The discrete-time dynamical systems approach to the analysis of the model assumed in this work   leads to consideration of a singular perturbation problem for which the unperturbed invariant manifold  is nonsmooth. Hence the existing theory for normally hyperbolic invariant manifolds (\cite{jones} and references therein) cannot be directly applied. The main result presented herein  is the proof of the persistence of an invariant manifold under singular perturbation and its consequences for understanding model behavior and  bifurcation scenarios. The technique used in the proof is Hadamard's graph transform method.

The components of the model, including an albedo function designed for the Neoproterozoic ice ages under consideration and following  that presented in \cite{abb}, are introduced in Section 2. In Section 3 the piecewise-defined model difference equations are described. While the system of equations is everywhere continuous, the system is shown to be nonsmooth on a  certain hyperplane in phase space. 

The model is placed within the context of  singular perturbation problems, with singular parameter $\eps$ representing the time constant for movement of the ice sheet, in Section 4. In addition, the unperturbed invariant manifold is shown to be Lipschitz continuous  and the function spaces to be used for the graph transforms are   introduced in Section 4. In Section 5 the existence of an attracting invariant manifold is proved to exist, for sufficiently small $\eps$, and in Section 6 we present the model dynamics and we discuss local smooth and nonsmooth bifurcations. We provide concluding remarks in Section 7.

\section{The temperature--ice line model}

Let $T_n(y) \ (^\circ$C) denote the average annual surface temperature at $y=\sin\theta$, where $\theta$ is the latitude, at year $n$. The variable $y$, referred to as ``latitude" in all that follows, is convenient in that a latitudinal band at $y$ has area proportional to $dy$. The energy balance equation for temperature  has three terms, corresponding to absorbed insolation, outgoing longwave radiation, and meridional heat transport, respectively. A diffusion approach for the heat transport term is used here, as is the case, for example, in \cite{cal}, \cite{north2}, \cite{north}, \cite{north84} and \cite{northetal}. Meridional heat transport refers to the heat flux across latitude circles associated with   atmospheric and oceanic currents.

Consider the difference equation
\begin{equation}\label{Tnsetup}
R\cdot\frac{T_{n+1}(y)-T_n(y)}{(n+1)-n}=Qs(y)(1-\al(y,\eta))-(A+BT_n(y))+D\frac{\partial}{\partial y} (1-y^2)\frac{\partial T_n}{\partial y}.
\end{equation}
The model assumes a symmetry across the equator, so we take $y\in[0,1]$.  
The units on each side of \eqref{Tnsetup} are Wm$^{-2}$. \ $R$ is the heat capacity of the Earth's surface, with units Wyr$(^\circ$C m$^2)^{-1}$, while $Q$   (Wm$^{-2}$) is the average annual insolation for the entire Earth. The function
\begin{equation}\notag
s(y)=\frac{2}{\pi^2}\int^{2\pi}_0 \sqrt{1-(\sqrt{1-y^2} \, \sin\beta\cos\gamma-y\cos\beta)^2} \, d\gamma,
\end{equation}
in which $\beta$ denotes the tilt (or {\em obliquity}) of the Earth's spin axis, serves to incorporate insolation as a function of latitude \cite{dickclar}. For the  planet Earth, $s(y)$ monotonically decreases from a maximum $s(0)$ at the equator to a minimum $s(1)$ at the North Pole. 

The model assumes ice exists everywhere above the edge of the ice sheet, denoted $y=\eta$ and called the {\em ice line}, with no ice at latitudes equatorward of $y=\eta$. The function $\al(y,\eta)$ represents the surface albedo, which depends upon the position of the ice line. Hence $Qs(y)(1-\al(y,\eta))$ models the incoming solar radiation absorbed at the surface at latitude $y$.

The $(A+BT_n)$-term represents the outgoing longwave radiation, with $A  $ (Wm$^{-2}$) and $B$ (W($^\circ$C m$^2)^{-1}$) estimated by satellite data \cite{graves}.  The diffusive meridional heat transport term is assumed to be proportional to the surface temperature gradient, and  
\begin{equation}\label{diffuse}
D\frac{\partial}{\partial y} (1-y^2)\frac{\partial T}{\partial y}
\end{equation} is the form the spherical diffusion operator takes when assuming no radial or longitudinal dependence. The parameter $D>0 \, $ (W($^\circ$C m$^2)^{-1}$) is a  diffusion coefficient. The diffusive approach comes with the boundary conditions that the gradiant of the temperature profile $T_n(y)$ equals zero at the equator and at the North Pole.

We note that  Budyko  \cite{bud} and Sellers \cite{sell} each utilized a discrete time approach in their seminal 1969 papers. We also note Budyko used an alternative formulation of the heat transport term, namely,
\begin{equation}\label{mean}
-C\left(T(y)-\tx{\int^1_0 } T(y) dy\right), \ C>0.
\end{equation}
In this approach it is the temperature at latitude $y$ relative to the global average surface temperature that governs the meridional heat transport.

Equation \eqref{Tnsetup} leads to consideration of the temperature evolution equation
\begin{align}\label{Tdiff}
T_{n+1}(y)=T_n(y)+\frac{1}{R}\left(Qs(y)(1-\al(y,\eta))-(A+BT_n(y))+D\frac{\partial}{\partial y} (1-y^2)\frac{\partial T_n}{\partial y}\right).
\end{align}
Budyko, Sellers and others working with energy balance climate models focused on the  equilibrium temperature profiles $T^*(y)$
of model equations such as \eqref{Tdiff}, and how these solutions varied with a parameter (typically $Q$). In particular, there was no consideration of an ice line allowed to move with changes in temperature. Not until E. Widiasih's paper \cite{esther} was the evolution of surface temperature  coupled to a dynamic ice line, thereby introducing  a dynamical systems approach to the study of the model. Although Widiasih used   heat transport term \eqref{mean} and worked in the infinite dimensional setting (with parameters aligned with the present climate), the methods used here are similar in spirit to those  in \cite{esther}. 

We consider the coupled system
\begin{subequations}\label{eq:budwid}
\begin{align}
T_{n+1}(y)& =T_n(y)+\frac{1}{R}\left(Qs(y)(1-\al(y,\eta_n))-(A+BT_n(y))+D\frac{\partial}{\partial y} (1-y^2)\frac{\partial T_n}{\partial y}\right) \label{eq:budwidA}\\  
\eta_{n+1}&= \eta_n+\eps(T_n(\eta_n)-T_c), \label{eq:budwidB}
\end{align}
\end{subequations}
where $\eps>0$ is a small  parameter and $T_c$ is a {\em critical temperature}. Note that $T_n(\eta_n)>T_c$ implies $\eta_{n+1}>\eta_n$, so that the ice line moves poleward. If $T_n(\eta_n)<T_c$, then $\eta_{n+1}<\eta$ and the ice line descends equatorward. Equilibria of system \eqref{eq:budwid} are pairs $(T^*(y),\eta^*)$ for which  equation \eqref{Tnsetup} equals zero, with the temperature at the ice line additionally satisfying $T^*(\eta^*)=T_c$.
The ice line quation \eqref{eq:budwidB} was introduced in \cite{esther}.

\subsection{The albedo function}

The albedo function used here is motivated by \cite{abb} , in which the extensive glacial episodes of the Neoproterozoic  Era were modeled. Utilizing an idealized general circulation model parametrized for these great ice ages, Abbot et al found that evaporation exceeded precipitation in a latitude band  situated within the tropics. This leads to the conceptual model assumption that ice forming below a fixed latitude $y=\rho$ would not have a snow cover, resulting in an albedo for this ``bare" ice lying between the high albedo of snow covered ice and the low albedo of ice free surface. There is then an increase in absorbed solar radiation through the bare ice, relative to that passing through the snow covered ice. 
The resulting increase in local surface temperature  provides a possible mechanism by which  the advance of the ice sheet might be arrested. 

A stable Jormungand equilibrium solution, with the ice line resting in tropical latitudes, was found by Abbot et al when running the large computer model used in \cite{abb}. Working in the infinite dimensional setting and using heat transport term \eqref{mean}, a dynamically stable Jormungand equilibrium solution was shown to exist in \cite{mewid} when using a smooth version of albedo function \eqref{Jalb}.

Incorporating the lower albedo of bare ice into \eqref{eq:budwid}, we fix a latitude $y=\rho$ as in \cite{abb}  and we define the albedo function as follows:
\begin{equation} \label{Jalb}
\al(y,\eta)=
\begin{cases}
\al^-(y,\eta), & \eta<\rho \\
\al^+(y,\eta), & \rho\leq\eta,\\
\end{cases}
\end{equation}
where
\begin{equation} \label{andb}
\begin{array}{cc}
\al^-(y,\eta)=
\begin{cases}
\al_1, & y<\eta \\
\frac{1}{2}(\al_1+\al_i),  & y=\eta \\
\al_i, & \eta<y<\rho\\
\frac{1}{2}(\al_i+\al_2),  & y=\rho \\
\al_2, & \rho<y,
\end{cases}  & \hspace*{.1in} \mbox{and}  \hspace*{.15in}
\al^+(y,\eta)=
\begin{cases}
\al_1, & y<\eta \\
\frac{1}{2}(\al_1+\al_2),  & y=\eta \\
\al_2, & y>\eta.
\end{cases}
\end{array}
\end{equation}
Here,  $\al_1<\al_i<\al_2$ represent the albedos of ice free surface, bare ice, and snow covered ice, respectively.
The use of the superscript $^+$ will indicate $\rho\leq\eta$, while the superscript $^-$ will indicate $\eta<\rho$,  in all that follows.

\section{The  model difference equations}

The use of Legendre polynomial expansions to analyze the  temperature equation \eqref{Tdiff} is quite common  (\cite{hans}, \cite{north2}, \cite{north}, \cite{north84}, \cite{northetal}), due to the fact the Legendre polynomials $p_n(y)$ are  eigenfunctions of the spherical diffusion operator:
\begin{equation}\notag
\frac{d}{dy}(1-y^2)\frac{d}{dy}p_n(y)=-n(n+1)p_n(y), \ n=0, 1, ... \, .
\end{equation}
We incorporate finite Legendre expansion  approximations in the coupled temperature-ice line model \eqref{eq:budwid}. Recall $T_n(y), s(y)$ and $\al(y,\eta)$ are each even functions of $y$, given the assumed symmetry of the model about the equator.
We thus write
\begin{equation}\label{LegT}
T_n(y)=\sum^N_{i=0} x_{n, 2i}p_{2i}(y),
\end{equation}
with $x_{n,2i}$ undetermined coefficients, $i=0, ... , N$. As the even Legendre polynomials each have zero gradient at the equator and at the North Pole \cite{north}, $T_n(y)$ satisfies the prescribed boundary conditions for any $N$.

We also write
\begin{equation}\label{ess}
s(y)= \sum^N_{i=0} s_{2i} \, p_{2i}(y),
\end{equation}
where
\begin{equation}\notag
s_{2i}=(4i+1) \int^1_0s(y)p_{2i}(y)dy.
\end{equation}
Finally, we express
\begin{equation}\label{essalpha}
s(y)\al^\pm(y,\eta)=\sum^N_{i=0} a^\pm_{2i}(\eta)p_{2i}(y) dy,
\end{equation}
where 
\begin{align}\notag
a^+_{2i}(\eta)&=(4i+1) \int^1_0s(y)\al^+(y,\eta)p_{2i}(y)dy=\al_2s_{2i}-(4i+1)(\al_2-\al_1) \int^\eta_0s(y)p_{2i}(y)dy, 
\end{align}
and
\begin{align}\notag
a^-_{2i}(\eta)&=(4i+1) \int^1_0s(y)\al^-(y,\eta)p_{2i}(y)dy\\\notag
&=\al_2s_{2i}-(4i+1)\left( (\al_2-\al_i) \int^\rho_\eta s(y)p_{2i}(y)dy+(\al_2-\al_1) \int^\eta_0s(y)p_{2i}(y)dy \right).\notag
\end{align}
Again, $a^+_{2i}(\eta)$ corresponds to $\rho\leq \eta$, while  $a^-_{2i}(\eta)$ corresponds to $\eta<\rho$, via definition \eqref{Jalb}.

Substitution of \eqref{LegT}, \eqref{ess} and \eqref{essalpha} into equation  \eqref{eq:budwidA} yields
\begin{align}\label{subst}
\tx{\sum^N_{i=0}} x_{n+1,2i} \, p_{2i}(y)&=\tx{\sum^N_{i=0}} x_{n,2i} \, p_{2i}(y)+\frac{1}{R}\left[Q\tx{\sum^N_{i=0}}(s_{2i}-a^\pm_{2i}(\eta)) \, p_{2i}(y)\right.\\\notag
&-\left(A+B \tx{\sum^N_{i=0}}x_{n,2i} \, p_{2i}(y)\right)
\left.-D \, \tx{\sum^N_{i=0}} 2i(2i+1)x_{n,2i} \, p_{2i}(y)\right].
\end{align}
Equating the coefficients of $p_{2i}(y)$ in \eqref{subst}, we have
\begin{align}\notag
x_{n+1,0}&=x_{n,0}+\tx{\frac{1}{R}}\left[Q(s_0-a^\pm_0(\eta))-A-Bx_{n,0}\right],\\\notag
x_{n+1,2i}&=x_{n,2i}+\tx{\frac{1}{R}}\left[Q(s_{2i}-a^\pm_{2i}(\eta))-(B+2i(2i+1)D)x_{n,2i}\right], \ i=1, ... , N.
\end{align}

Adding the ice line equation \eqref{eq:budwidB}
and doing a bit of rearranging,   one finds the $N+2$ piecewise-defined equations
\begin{align}\label{mysystem}
x_{n+1,2i}&= x_{n,2i}-\textstyle{\frac{(B+2i(2i+1)D)}{R}}\left(x_{n,2i}-f^\pm_{2i}(\eta_n)\right), \ i=0, ... , N,\\\notag
\eta_{n+1}&=\eta_n+\eps\left(\textstyle{\sum^N_{i=0}} \, x_{n,2i} \, p_{2i}(\eta_n)-T_c\right),\notag
\end{align}
where
\begin{equation}\notag
f^\pm_0(\eta)=\textstyle{\frac{1}{B}}(Q(s_0-a^\pm_0(\eta))-A),  \ \ f^\pm_{2i}(\eta)=\textstyle{\frac{1}{B+2i(2i+1)D}} \, Q(s_{2i}-a^\pm_{2i}(\eta)), \ i\geq 1.
\end{equation}
(Recall $f^+_{2i}(\eta)$ corresponds to $\rho\leq\eta$, and $f^-_{2i}(\eta)$ corresponds to $ \eta<\rho$.) 
 We  note each $f^\pm_{2i}(\eta)$ is a polynomial of degree $2N+2i+1$.

The state space for \eqref{mysystem} is $ \R^{N+1}\times [0,1]$. We observe that $a^+_{2i}(\rho)=a^-_{2i}(\rho),$ so that $f^+_{2i}(\rho)=f^-_{2i}(\rho), \, i=0, ... , N$, and hence \eqref{mysystem} is everywhere continuous. Notably, however, system \eqref{mysystem} is not smooth at any point in the set 
\begin{equation}\label{Sigma}
\Sigma=\{ (x,\eta) : x\in\R^{N+1}, \ \eta=\rho\};
\end{equation}
  for example, $\frac{\partial f^-_0}{\partial \eta}=\frac{Q}{B}(\al_i-\al_1)s(\eta)$ and $\frac{\partial f^+_0}{\partial \eta}=\frac{Q}{B}(\al_2-\al_1)s(\eta)$ do not agree on $\Sigma.$ This fact provides for the lack of smoothness of the invariant manifolds to be described below.

For convenience with the ensuing analysis, we extend system \eqref{mysystem} to $\R^{N+1}\times\R$  as follows. (For the persistence of smooth, noncompact normally hyperbolic invariant manifolds  of bounded geometry, see \cite{jaap}.) 
Extend  expression \eqref{LegT} to $y\in\R$ by setting $T_n(y)=\sum^N_{i=0}x_{n,2i}p_{2i}(0)$ for $y<0$, and by setting $T_n(y)=\sum^N_{i=0}x_{n,2i}p_{2i}(1)$ for $y>1$. That is, for all $y<0$ we set $T_n(y)=T_n(0)$, while  for all $y>1$ we set $T_n(y)=T_n(1)$.
 Note albedo function \eqref{andb} is well defined for $y,\eta\in\R$. 
For $i=0, ... , N$, set
\begin{equation}\label{extendf}
\begin{array}{cc}
f_{2i}(\eta)=
\begin{cases}
f^-_{2i}(0), & \eta<0\\
f^-_{2i}(\eta), & 0\leq \eta<\rho\\
f^+_{2i}(\eta), & \rho\leq \eta\leq 1\\
f^+_{2i}(1), & 1<\eta,
\end{cases} & \hspace*{.2in} \mbox{and}  \hspace*{.2in}
q_{2i}(\eta)=
\begin{cases}
p_{2i}(0), & \eta<0 \\
p_{2i}(\eta), & 0\leq\eta\leq 1\\
p_{2i}(1), & 1<\eta.
\end{cases}
\end{array}
\end{equation}
These choices of $f_{2i}(\eta)$ and $q_{2i}(\eta)$ ensure that
\begin{align}\label{extsystem}
x_{n+1,2i}&= x_{n,2i}-\textstyle{\frac{(B+2i(2i+1)D)}{R}}\left(x_{n,2i}-f_{2i}(\eta_n)\right), \ i=0, ... , N\\\notag
\eta_{n+1}&=\eta_n+\eps\left(\textstyle{\sum^N_{i=0}} \, x_{n,2i} \, q_{2i}(\eta_n)-T_c\right)\notag
\end{align}
is a continuous system of difference equations for $(x,\eta)\in\R^{N+1}\times\R$. We note this system is nonsmooth at any point $(x,\eta)$ with $\eta=0,\rho$ or 1. 

For $v \in \R^{N+1}$, let $(v)_i$ denote the $i$th coordinate  of $v$, where we take $i=0, 1, ... , N$. For ease of notation, set 
\begin{equation}\notag
\gamma_{i}=\frac{B+2i(2i+1)D}{R},  \ i=0, ... , N.
\end{equation}
Consider the functions
\begin{equation}\label{F}
F :\R^{N+1}\times\R\rw \R^{N+1}, \ (F(x,\eta))_{i}=-\gamma_{i}((x)_{i}-f_{2i}(\eta)), \ i=0, ... , N, 
\end{equation}
and
\begin{equation}\label{Gee}
G:\R^{N+1}\times\R\rw\R,  \ G(x,\eta)=\eps\left( \tx{\sum^N_{i=0} } \, (x)_i \, q_{2i}(\eta)-T_c\right).
\end{equation}
Letting $X_n=(x_{n,0}, ... , x_{n,2N}),$ system \eqref{extsystem} is equivalent to 
\begin{equation}\label{EffGee}
(X_{n+1},\eta_{n+1})=(X_n+F(X_n,\eta_n), \ \eta_n+G(X_n,\eta_n)).
\end{equation}
We therefore are interested in understanding the behavior of orbits  for the mapping 
\begin{equation}\label{H}
H: \R^{N+1}\times\R\rw \R^{N+1}\times\R, \ H(x,\eta)=(x+F(x,\eta), \eta+G(x,\eta)). 
\end{equation}

\begin{figure}[t]
\includegraphics[width=6.1in,trim = 1.3in 7.5in .5in  1.1in, clip]{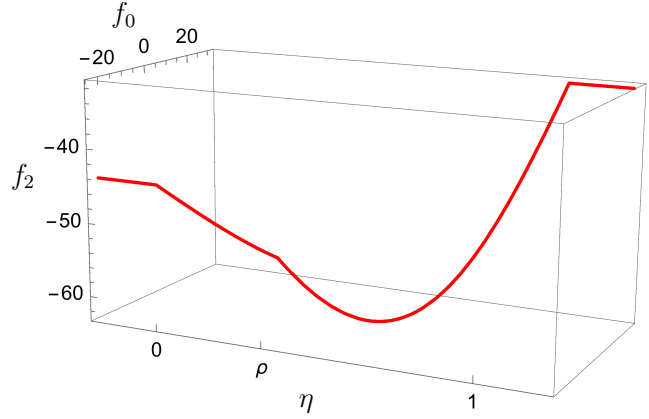}
\begin{center}
\caption{The invariant manifold of fixed points $M_0=\{(h^0(\eta),\eta)=(f_0(\eta), f_2(\eta),\eta) : \eta\in\R\}$ for $N=1$.
}
 \end{center}
\end{figure}

\section{The case $\boldsymbol{\eps=0}$}

\subsection{The invariant manifold is Lipschitz continuous} 
The position of the ice line remains fixed for all time when $\eps=0$. In this case \eqref{extsystem} admits  a manifold $M_0=\mbox{graph}(h^0)$ of fixed points, where
\begin{equation}\label{h0}
h^0 :\R\rw \R^{N+1}, \ \eta\mapsto (f_0(\eta), f_2(\eta), ... , f_{2N}(\eta))
\end{equation}
(see Figure 1). We assume throughout that $N$ is (maximally) chosen so that
\begin{equation}\label{Nbound}
\left|1-\gamma_N\right|= \left|1-\textstyle{\frac{B+2N(2N+1)D}{R}}\right|\leq 1-\textstyle{\frac{B}{R}}=1-\gamma_0<1.
\end{equation}
Note \eqref{Nbound} implies $|1-\gamma_i|\leq 1-\gamma_0<1$ for $i=0, ... , N$. 
With $R, B$ and $D$ values as in Table 1,  $N=5$ satisfies requirement \eqref{Nbound}, ensuring that expansion \eqref{LegT}  provides a good approximation to the temperature distribution function $T_n(y)$. See Figure 2, for example, in which the equilibrium temperature profiles
\begin{equation}\notag
T^*(y)=\tx{\sum^N_{i=0}} f_{2i}(\eta)p_{2i}(y)
\end{equation}
are plotted when $\eps=0$ for $\eta=0.3$ and $N=4, 5$.

Assumption \eqref{Nbound} guarantees that for $\eps=0$, \, $M_0$ is a globally attracting curve of fixed points. As noted above, the existing theory for the persistence of normally hyperbolic  invariant manifolds under perturbation cannot be applied due to the lack of smoothness of $M_0$. We will use Hadamard's graph transform technique to prove the existence of an attracting invariant manifold $M_\eps$ that is within $O(\eps)$ of $M_0$, for sufficiently small $\eps$. This will allow for the analysis of the behavior of orbits of the function $H(x,\eta)$ \eqref{H}, and provide for an understanding of  bifurcations for system \eqref{extsystem}.

We first note that, while nonsmooth, $M_0$ is the graph of a Lipschitz continuous function.

\begin{table}[t]
\small
\centering
\begin{tabular}{|ccc|ccc|}\hline
Parameter & Value & Units & Parameter & Value & Units\\ \hline
$R$ & 20   & Wyr$(^\circ$C m$^2)^{-1}$  & $\al_i$ & 0.4 & dimensionless\\
$Q$  & 321  & Wm$^{-2}$ & $\al_2$ & 0.8 & dimensionless\\ 
$A$ & 164   & Wm$^{-2}$ & $s_0$ & 1 & dimensionless\\ 
$B$ & 1.9  & W$(^\circ$C m$^2)^{-1}$ & $s_2$ & -0.477131 & dimensionless\\
$D$ & 0.25   & W$(^\circ$C m$^2)^{-1}$  & $s_4$ & -0.045029 & dimensionless\\
$\beta$ & 23.4   & degrees & $s_6$ &  0.007937 & dimensionless\\
$T_c$ & 0    & $^\circ$C & $s_8$ & 0.013859 & dimensionless\\
 $\al_1$ & 0.30 & dimensionless & $s_{10}$ & 0.008663 & dimensionless\\ \hline
\end{tabular}
\caption{Values of parameters and constants.}

\vspace{.1in}
\end{table}

\vspace{0.25in}
\noindent
{\em Proposition 1}. The function $h^0(\eta)$ is Lipschitz continuous.

\vspace{0.15in}
\noindent
{\em Proof.} \ Let $h^0_\pm(\eta)=(f^\pm_0(\eta), ... , f^\pm_{2N}(\eta))$. For $\rho\leq \eta\leq1,$ we have
\begin{align}\notag
\|(h^0_+)^\pr(\eta)\|&=\left[ \tx{\sum^N_{i=0}} \left( \frac{(4i+1)Q}{(B+2i(2i+1)D)}(\al_2-\al_1)s(\eta)q_{2i}(\eta)\right)^2\right]^{1/2}\notag
\end{align}
\begin{align}\notag
&\leq \left[\tx{\sum^N_{i=0}}\left(\frac{(4N+1)Q}{B}(\al_2-\al_1)s(0)\right)^2\right]^{1/2}\\\notag
&=\frac{(4N+1)Q}{B}(\al_2-\al_1)s(0)\sqrt{N+1}. 
\end{align}
Note we used $s(0)$ and 1 as upper bounds for $s(\eta)$ and $q_{2i}(\eta)$, respectively, on $[0,1]$.

\begin{figure}[t]
\includegraphics[width=6.1in,trim = 1.3in 7.75in .5in  1.25in, clip]{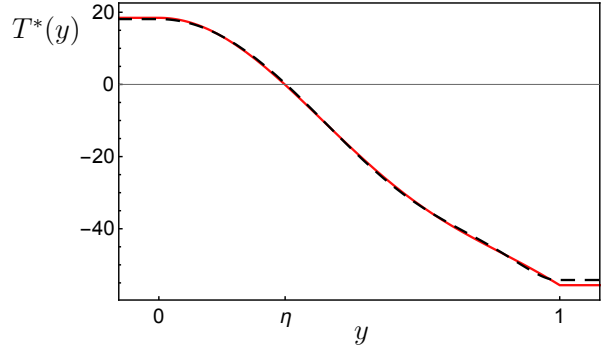}
\begin{center}
\caption{Plots of the equilibrium temperature profile $T^*(y)$ when $\eps=0$ for $\eta=0.3$. {\em Red}: $N=4.$ \ {\em Dashed}: $N=5$.
}
 \end{center}
\end{figure}

For $0\leq  \eta\leq\rho,$ a similar computation yields
\begin{align}\notag
\|(h^0_-)^\pr(\eta)\|&\leq \frac{(4N+1)Q}{B}(\al_i-\al_1)s(0)\sqrt{N+1}. 
\end{align}
Suppose $0< \eta_1<\rho<\eta_2<1$. As $h^0$ is not smooth at points for which  $\eta=\rho$, the derivative cannot be used to bound Lip$(h^0$)   in this case.
Alternatively, and using $h^0_-(\rho)=h^0_+(\rho)$, we have
\begin{align}\notag
 \frac{\|h^0_+(\eta_2)-h^0_-(\eta_1)\|}{|\eta_2-\eta_1|}&\leq  \frac{\|h^0_+(\eta_2)-h^0_+(\rho)\|}{|\eta_2-\eta_1|}+
  \frac{\|h^0_-(\rho)-h^0_-(\eta_1)\|}{|\eta_2-\eta_1|}\\\notag
  &\leq  \frac{\|h^0_+(\eta_2)-h^0_+(\rho)\|}{|\eta_2-\rho|}+
  \frac{\|h^0_-(\rho)-h^0_-(\eta_1)\|}{|\rho-\eta_1|}\\\notag
  &\leq \mbox{Lip}_{[\rho,1]}(h^0_+)+\mbox{Lip}_{[0,\rho]}(h^0_-)\\\notag
  &\leq \frac{(4N+1)Q}{B}s(0)\sqrt{N+1} \, (\al_2-\al_1+\al_i-\al_1).
\end{align}
Let
 \begin{equation}\label{L0}
 L_0=\frac{(4N+1)Qs(0)}{B}\sqrt{N+1} \, (\al_2+\al_i-2\al_1).
 \end{equation}
One easily checks that $L_0$ serves as an upper bound for $\|h^0(\eta_2)-h^0(\eta_1)\|/|\eta_2-\eta_1|$ if $\eta_1<0$ or $1<\eta_2$ as well, given the simplicity of the extension \eqref{extsystem} of \eqref{mysystem}. Hence,
 Lip($h^0)\leq L_0. $ \ $_\square$

\subsection{Defining the function space ${\cal B}_{\bs{L}}$}

We introduce various quantities that will appear in ensuing estimates, as well as the function space to be used for the graph transform. Recall that $L_0 $ \eqref{L0} is a Lipschitz constant for $h^0(\eta)$. Let   
\begin{equation}\label{dee}
M=\|h^0\|_\infty, \  \ d\geq  1+2N(2N+1)\frac{D}{B}, 
\end{equation}
and set 
\begin{equation}\label{ell}
L=\max\{dL_0, dM\}. 
\end{equation}
We note $d\geq 1$ and,  given the parameters in Table 1, $L_0\geq 80$, so that $L\geq 80.$ 
Let
\begin{equation}\notag
K_i=\mbox{Lip}(q_{2i}), \ i=0, ... , N, \ \mbox{ and } \ K=\tx{\sum^N_{i=0}} \, K_i.
\end{equation}
We will make use of various inequalities involving the expressions above, including the following.

\vspace{0.25in}
\noindent
{\em Proposition 2}. (a) $-\gamma_0+\frac{\gamma_N}{d}\leq 0$.

\vspace{0.05in}
\noindent
(b) $ \gamma_0[L((1+\gamma_N)(N+1)+\gamma_0K)]^{-1}  \leq  [L(N+1+K)]^{-1}\leq  \gamma_0(N+1)^{-1}$.

\vspace{0.1in}
\noindent
(c) If \ $0<\eps< \gamma_0[L((1+\gamma_N)(N+1)+\gamma_0K)]^{-1}$, then 
\begin{equation}\notag
c= 1-\gamma_0+L(1-\gamma_0+\gamma_N)\tx{\frac{\eps(N+1)}{1-\eps L(N+1+K)}} \in[0,1).
\end{equation}

\vspace{0.1in}
\noindent
{\em Proof.} \ (a) Note $-\gamma_0+\frac{\gamma_N}{d}\leq 0$ is equivalent to $d\geq \frac{\gamma_N}{\gamma_0}=1+2N(2N+1)\frac{D}{B}$, which is assumed in \eqref{dee}. 

\vspace{0.1in}
\noindent
(b) Note $\frac{1+\gamma_N}{\gamma_0}\geq \frac{1}{\gamma_0}=\frac{R}{B}\geq 1$. Hence
\begin{equation}\notag
\frac{ \gamma_0}{L((1+\gamma_N)(N+1)+\gamma_0K)}=\frac{ 1}{L\left((\frac{1+\gamma_N}{\gamma_0})(N+1)+ K\right)}\leq \frac{1}{L(N+1+K) }.
\end{equation}
We also have
\begin{equation}\notag
\frac{1}{L(N+1+K) }\leq \frac{1}{L(N+1) }\leq \frac{\gamma_0}{N+1 },
\end{equation}
with the latter inequality following from $\frac{1}{L}\leq \frac{1}{80}\leq \frac{B}{R}=\gamma_0$.

\vspace{0.1in}
\noindent
(c)  Note  the hypothesis of part (c), together with part (b), imply $c>0$. Showing $c<1$ is equivalent to showing 
\begin{equation}\notag
L(1-\gamma_0+\gamma_N)\eps (N+1)<\gamma_0(1-\eps L(N+1+K)),
\end{equation}
or
\begin{equation}\notag
\eps\left[ L(N+1)(1-\gamma_0+\gamma_N)+\gamma_0L(N+1+K)\right]<\gamma_0,
\end{equation}
an expression equivalent to $\eps< \gamma_0[L((1+\gamma_N)(N+1)+\gamma_0K)]^{-1}$. \ $_\square$

\vspace{0.15in}
We will show the desired invariant manifold $M_\eps$ is the graph of a function $g^*(\eta)$, for $\eps>0$ and sufficiently small. To that end, we utilize the function spaces
\begin{equation}\notag
{\cal B}=\{g:\R\rw \R^{N+1} \  | \  g \mbox{ is continuous and }\|g\|_\infty<\infty\},
\end{equation}
and 
\begin{equation}\notag
{\cal B}_L=\{g\in {\cal B} \ | \  \|g\|_\infty\leq L \mbox{ and } \mbox{Lip}(g)\leq L\}.
\end{equation}
We turn to the main result in the following section.

\section{The invariant manifold $\bs{M_\eps}$}

In this section we prove the  following result.

\vspace{0.15in}
\noindent
{\em Theorem 1.} \ For $\eps>0$ and sufficiently small, there exists $g^*\in {\cal B}_L$ such that $M_\eps=\mbox{graph}(g^*)$ is a locally attracting invariant manifold for mapping \eqref{H}. In addition, $M_\eps$ lies within $O(\eps)$ of $M_0$.

\vspace{0.15in}
That  Theorem 1 holds true follows from the following propositions.

\vspace{0.2in}
\noindent
\noindent
{\em Proposition 3}. Assume $0<\eps< (L(N+1+K))^{-1}$. Let $g\in {\cal B}_L$. Then for all $\eta\in\R$, there exists $\beta=\beta(\eta,g)\in\R$ such that
\begin{equation}\label{preimage}
\eta=\beta+\eps\left(\tx{\sum^N_{i=0}} (g(\beta))_{i}q_{2i}(\beta)-T_c  \right).
\end{equation}

\vspace{0.05in}
\noindent
{\em Proof.} \ Let $g\in {\cal B}_L$. \ Consider the mapping
\begin{equation}\label{Phi}
\Phi:{\cal B}\rw {\cal B}, \ (\Phi b)(w)=\eps\left(T_c-\tx{\sum^N_{i=0}}(g(w+b(w)))_{i} \, q_{2i}(w+b(w))\right).
\end{equation}
Let $b_1, b_2\in {\cal B}_L$, let $w\in\R$ and, for ease of notation, let $u=w+b_2(w), \ v=w+b_1(w)$. We have
\begin{equation}\notag
\left|(\Phi b_1)(w)-(\Phi b_2)(w)\right|=\eps\left|\tx{\sum^N_{i=0}}[(g(u))_{i}q_{2i}(u)-(g(v))_{i}q_{2i}(v)]\right|.
\end{equation}
Using the fact $g\in {\cal B}_L$ and $\|q_{2i}\|_\infty=1$, we have, for $i=0, ... , N$,
\begin{align}\notag
|(g(u))_{i}q_{2i}(u) - (g(v))_{i} &q_{2i}(v)| \\\notag
&\leq |(g(u))_{i}q_{2i}(u)- (g(v))_{i}q_{2i}(u)|+|(g(v))_{i}q_{2i}(u)-(g(v))_{i}q_{2i}(v)|\\\notag
&\leq \mbox{Lip}((g)_{i})|u-v|+L K_i |u-v|\\\notag
&\leq L(1+K_i)|u-v|\\\notag
&=L(1+K_i)|b_2(w)-b_1(w)|\\\notag
&\leq L(1+K_i)\|b_2-b_1\|_\infty.
\end{align}
Since $w$ was arbitrary, we have
\begin{align}\notag
\|\Phi b_1-\Phi b_2\|_\infty\leq \eps \,  \tx{\sum^N_{i=0}}L(1+K_i)\|b_1-b_2\|_\infty=\eps L(N+1+K)\|b_1-b_2\|_\infty.
\end{align}
Give our assumed bound on $\eps$, we have shown $\Phi$ is a contraction map on the complete space (${\cal B},\| \cdot \|_\infty)$, implying $\Phi$ has a unique fixed point $b^*\in{\cal B}$ satisfying, for all $\eta\in\R$,
\begin{equation}\notag
b^*(\eta)=\eps\left(T_c-\tx{\sum^N_{i=0}}(g(\eta+b^*(\eta)))_{i} \, q_{2i}(\eta+b^*(\eta))\right)
\end{equation}
Letting $\beta=\eta+b^*(\eta)$ completes the proof.  \ $_\square$

\vspace{0.15in}
\noindent
{\em Remark}.   We note  $\beta$ in Proposition 3 depends upon  $g\in {\cal B}_L$, given that  the function $b^*$ in the proof of Proposition 3 depends on $g$. 
As in \cite{esther}, we call $\beta$ the {\em preimage} of the ice line $\eta$.

\vspace{0.15in}
The   next step is to produce a graph that is invariant under the  function $H(x,\eta)$.
 To that end, consider the mapping 
\begin{equation}\label{gamma}
\Gamma :{\cal B}_L\rw {\cal B}, \ (\Gamma g)(\eta)=g(\beta)+F(g(\beta),\beta),
\end{equation}
where $\beta$ is the preimage of  the ice line $\eta$ whose existence is guaranteed by  Proposition 3. We will show that $\Gamma$ fixes a function $g^*\in {\cal B}_L$, and that the graph of $g^*$ provides an  invariant manifold $M_\eps$ for $H(x,\eta)$, provided $\eps$ is sufficiently small.  

\vspace{0.25in}
\noindent
{\em Proposition 4}. Assume $0<\eps< \gamma_0[L((1+\gamma_N)(N+1)+\gamma_0K)]^{-1}$. There exists $c\in[0,1)$ such that for all $g,h\in {\cal B}_L$, \ 
$\|\Gamma g-\Gamma h\|_\infty\leq c\|g-h\|_\infty$.

\vspace{0.15in}
\noindent
{\em Proof.} \ Let $g, h\in {\cal B}_L$ and let $\eta\in\R$. Via Proposition 3, there exist $\beta_1, \beta_2\in\R$ such that
\begin{equation}\label{twobetas}
\eta=\beta_1+\eps\left(\tx{\sum^N_{i=0}} (g(\beta_1))_i \, q_{2i}(\beta_1)-T_c\right)
=\beta_2+\eps\left(\tx{\sum^N_{i=0}} (h(\beta_2))_i \, q_{2i}(\beta_2)-T_c\right).
\end{equation}
By definition, $(\Gamma g)(\eta)-(\Gamma h)(\eta)=g(\beta_1)+F(g(\beta_1),\beta_1)-h(\beta_2)-F(h(\beta_2),\beta_2)$, the $i$th component of which has the form
\begin{align}\notag
&(g(\beta_1))_i-\gamma_i( (g(\beta_1))_i-f_{2i}(\beta_1))-(h(\beta_2))_i+\gamma_i( (h(\beta_2))_i-f_{2i}(\beta_2))\\
&=(1-\gamma_i)[(g(\beta_1))_i-(h(\beta_2))_i]+\gamma_i(f_{2i}(\beta_1)-f_{2i}(\beta_2)),  
\end{align}
$i=0, ... , N.$ 
Recalling assumption \eqref{Nbound}, we have
\begin{align}\label{step1}
\|(\Gamma g)(\eta)-(\Gamma h)(\eta)\|&\leq (1-\gamma_0)\|g(\beta_1)-h(\beta_2)\|+\gamma_N\|h^0(\beta_1)-h^0(\beta_2)\|\\\notag
&\leq (1-\gamma_0)\left(\|g(\beta_1)-g(\beta_2)\|+\|g(\beta_2)-h(\beta_2)\|\right)+\gamma_NL_0|\beta_1-\beta_2|\\\notag
&\leq (1-\gamma_0)\left( L|\beta_1-\beta_2|+\|g-h\|_\infty\right)+\gamma_N L|\beta_1-\beta_2|\\\notag
&=(1-\gamma_0)\|g-h\|_\infty+L(1-\gamma_0+\gamma_N)|\beta_1-\beta_2|.
\end{align}
Note that
\begin{equation}\label{beta1-beta2}
|\beta_1-\beta_2|=\eps \left| \tx{\sum^N_{i=0}}[(h(\beta_2))_iq_{2i}(\beta_2)-(g(\beta_1))_iq_{2i}(\beta_1)]\right|.
\end{equation}
Moreover  for $i=0, ... , N$,
\begin{align}\notag
|(h(\beta_2))_iq_{2i}(\beta_2&)-(g(\beta_1))_iq_{2i}(\beta_1)|\\\notag
&\leq|(h(\beta_2))_iq_{2i}(\beta_2)-(h(\beta_2))_iq_{2i}(\beta_1)|+|(h(\beta_2))_iq_{2i}(\beta_1)
-(g(\beta_1))_iq_{2i}(\beta_1)|\\\notag
&\leq \|h\|_\infty K_i|\beta_2-\beta_1|+\|h(\beta_2)-g(\beta_1)\|\\\notag
&\leq LK_i|\beta_2-\beta_1|+\|h(\beta_2)-h(\beta_1)\|+\|h(\beta_1)-g(\beta_1)\|\\\notag
&\leq LK_i|\beta_2-\beta_1|+L|\beta_1-\beta_2|+\|g-h\|_\infty\\\notag
&=L(K_i+1)|\beta_1-\beta_2|+\|g-h\|_\infty.
\end{align}
Returning to \eqref{beta1-beta2}, we have
\begin{align}\label{morebetas}
|\beta_1-\beta_2|&\leq \eps \tx{\sum^N_{i=0}}[L(K_i+1)|\beta_1-\beta_2|+\|g-h\|_\infty].
\end{align}
Invoking part (b) of Proposition 2, inequality
  \eqref{morebetas} simplifies to
\begin{equation}\label{evenmore}
|\beta_1-\beta_2|\leq  \frac{\eps(N+1)}{1-\eps L(N+1+K)}\|g-h\|_\infty.
\end{equation}
Substituting \eqref{evenmore} into \eqref{step1} yields
\begin{equation}\notag
\|(\Gamma g)(\eta)-(\Gamma h)(\eta)\|\leq \left(1-\gamma_0+L(1-\gamma_0+\gamma_N)\tx{\frac{\eps(N+1)}{1-\eps L(N+1+K)}}\right)\|g-h\|_\infty.
\end{equation}
Note
\begin{equation}\notag
c= 1-\gamma_0+L(1-\gamma_0+\gamma_N)\tx{\frac{\eps(N+1)}{1-\eps L(N+1+K)}} \in[0,1)
\end{equation} 
by Proposition 2, part (c). 
Hence $\|\Gamma g-\Gamma h\|_\infty\leq c \|g-h\|_\infty$, as desired. \ $_\square$

\vspace{0.15in}
We highlight the conditions imposed on $\eps$ in Proposition 4:
\begin{equation}\label{epsbound}
0<\eps< \gamma_0[L((1+\gamma_N)(N+1)+\gamma_0K)]^{-1}.
\end{equation}

\vspace{0.1in}
\noindent
{\em Proposition 5}. For $\eps$ satisfying \eqref{epsbound}, the mapping $\Gamma :{\cal B}_L\rw {\cal B}$ has a unique fixed point $g^*\in {\cal B}_L$.

\vspace{0.2in}
\noindent
{\em Proof.} \ We show that for all $g\in {\cal B}_L, \ \Gamma g\in  {\cal B}_L$, and then use the fact  ${\cal B}_L$ is a closed subset of the complete space $({\cal B}, \| \cdot \|_\infty)$. The conclusion then  follows by invoking  Proposition 4.

Let $g\in  {\cal B}_L, \ \eta\in\R$. Letting $\beta\in\R$ denote the preimage of $\eta$ for $g$, as in Proposition 3, we have $(\Gamma g)(\eta) =g(\beta)+F(g(\beta),\beta)$. Note 
\begin{equation}\label{boundGamma}
((\Gamma g)(\eta))_i= (1-\gamma_i)(g(\beta))_i+\gamma_if_{2i}(\beta), \ i=0, ... , N.
\end{equation}
It then follows from \eqref{boundGamma}  and assumption \eqref{Nbound}   that  
\begin{align}\notag
\|(\Gamma g)(\eta)\|&\leq  (1-\gamma_0)\|g(\beta)\|+\gamma_N\|h^0(\beta)\|\\\notag
&\leq (1-\gamma_0)\|g\|_\infty +\gamma_N\|h^0\|_\infty\\\notag
&\leq (1-\gamma_0)L+\gamma_NM\leq (1-\gamma_0+\tx{\frac{\gamma_N}{d}})L,
\end{align}
using the fact $L\geq dM$.
That $1-\gamma_0+\tx{\frac{\gamma_N}{d}}\leq 1$ follows from Proposition 2, part (a). Since $\eta$ was arbitrary, we have 
 $\|\Gamma g\|_\infty\leq L$. 

Using the fact the Lipschitz constant of the sum of two functions is the sum of the Lipschitz constants, together with \eqref{boundGamma} and \eqref{Nbound}, we have
\begin{align}\notag
\mbox{Lip}(\Gamma g)&\leq (1-\gamma_0)\mbox{Lip}(g)+\gamma_N\mbox{Lip}(h^0)\\\notag
&\leq (1-\gamma_0)L+\gamma_NL_0
\\\notag
&\leq (1-\gamma_0+\tx{\frac{\gamma_N}{d}})L\leq L,
\end{align}
using the fact $L\geq dL_0$. 
We have shown  $\Gamma g\in{\cal B}_L$. \ $_\square$

\vspace{0.25in}
\noindent
{\em Proposition 6}.  For $\eps$ satisfying \eqref{epsbound},  $M_\eps=\mbox{graph}(g^*)$ is invariant under the mapping $H$, where 
 $g^*\in{\cal B}_L$ is the unique fixed point of $ \, \Gamma$.

\vspace{0.15in}
\noindent
{\em Proof.} \ Let $\beta\in\R$. Let $\eta\in\R$ satisfy $\beta=\eta+b_{g^*}(\eta)$, where $b_{g^*}$ is the fixed point of the mapping $\Phi$, corresponding to $g^*$, introduced in the proof of Proposition 3. Then $\eta=\beta+G(g^*(\beta),\beta)$, and by Proposition 5
\begin{equation}\notag
(\Gamma g^*)(\eta)=g^*(\beta)+F(g^*(\beta),\beta)=g^*(\eta).
\end{equation}
Hence \ $H(g^*(\beta),\beta)=(g^*(\beta)+F(g^*(\beta),\beta), \ \beta+G(g^*(\beta),\beta))=(g^*(\eta),\eta). $ \ $_\square$

\vspace{0.25in}
\noindent
{\em Proposition 7}.  Assume  $\eps$ satisfies \eqref{epsbound} and let $M_\eps$ be as in Proposition 6. Then
\begin{itemize}
\item[(a)] $M_\eps$ is within $O(\eps)$ of $M_0$, and 
\item[(b)] For $x\in\R^{N+1}$ with $\|x\|\leq L$, and for $\beta\in\R$, \ $\|H^j(x,\beta)-M_\eps\|\rw 0$ as $j\rw\infty$, where $H^j$ denotes the $j$th iterate of $H$.
\end{itemize}

\vspace{0.15in}
\noindent
{\em Proof.} \ To prove (a), let $\eta\in\R$ and let $\beta=\eta+b_{g^*}(\eta)$, where
\begin{equation}\notag
 b_{g^*}(\eta)=\eps\left(T_c-\tx{\sum^N_{i=0}} \ (g^*(\beta))_i \, q_{2i}(\beta)\right),
\end{equation}
 as in Proposition 3. Recall
\begin{equation}\notag
g^*(\eta)=(\Gamma g^*)(\eta)=g^*(\beta)+F(g^*(\beta),\beta).
\end{equation}
By definition of $F$,
\begin{equation}\notag
(F(g^*(\beta),\beta))_i=-\gamma_i((g^*(\beta))_i-f_{2i}(\beta)).
\end{equation}
Once again invoking \eqref{Nbound}, we have
\begin{equation}\notag
\|g^*(\eta)-g^*(\beta)\|=\|F(g^*(\beta),\beta)\|\geq \gamma_0\|g^*(\beta)-h^0(\beta)\|,
\end{equation}
from which it follows

\begin{align}\notag
\|g^*(\beta)-h^0(\beta)\|&\leq \tx{\frac{1}{\gamma_0}}\|g^*(\eta)-g^*(\beta)\|=\tx{\frac{1}{\gamma_0}}\|g^*(\eta)-g^*(\eta+b_{g^*}(\eta))\|\\\notag
&\leq \tx{\frac{1}{\gamma_0}}L|b_{g^*}(\eta)|\leq \eps \tx{\frac{1}{\gamma_0}}L\left|T_c-\tx{\sum^N_{i=0}}(g^*(\beta))_i q_{2i}(\beta)\right|\\\notag
&\leq  \eps \tx{\frac{1}{\gamma_0}}L\left(|T_c|+\tx{\sum^N_{i=0}}L \right)=\eps \, \omega,
\end{align}
where $\omega$ is constant. Hence $M_\eps$ is within $O(\eps)$ of $M_0$.

For part (b), let $x\in \R^{N+1}$ with $\|x\|\leq L$. Let $\beta\in\R$, and set $p=(x,\beta)$. Let $g:\R\rw\R^{N+1}, \ g(\beta)=x$, so that $p=(g(\beta),\beta).$
Note $g\in {\cal B}_L$. Let $q=(g^*(\beta),\beta)$. Again as in the proof of Proposition 3, let $r,s\in\R$ satisfy $\beta=r+b_{g^*}(r)=s+b_g(s)$, that is,
\begin{equation}\notag
r=\beta+\eps\left(\tx{\sum^N_{i=0}}(g^*(\beta))_iq_{2i}(\beta)-T_c\right) \ \mbox{ and } \ s=\beta+\eps\left(\tx{\sum^N_{i=0}}(g(\beta))_iq_{2i}(\beta)-T_c\right).
\end{equation}
We have
\begin{align}\notag
|r-s|&=\eps\left| \tx{\sum^N_{i=0}}((g^*(\beta))_i-(g(\beta))_i)q_{2i}(\beta)\right|\\\notag
&\leq \eps \, \tx{\sum^N_{i=0}} |(g^*(\beta))_i-(g(\beta))_i|\\\notag
&\leq \eps (N+1)\|g^*(\beta)-g(\beta)\|=\eps (N+1)\|p-q\|.
\end{align}
Moreover, recalling  definitions \eqref{F} and \eqref{gamma}, one finds
\begin{equation}\notag
((\Gamma g^*)(r))_i-((\Gamma g)(s))_i=(1-\gamma_i)((g^*(\beta))_i-(g(\beta))_i),
\end{equation}
from which it follows
\begin{equation}\notag
\|(\Gamma g^*)(r)-(\Gamma g)(s)\|\leq (1-\gamma_0)\|g^*(\beta)-g(\beta)\|=(1-\gamma_0)\|p-q\|.
\end{equation}
We then have
\begin{align}\notag
\|H(p)-H(q)\|&=\|((\Gamma g^*)(r),r)-((\Gamma g)(s),s)\|\\\notag
&\leq \|(\Gamma g^*)(r)-(\Gamma g)(s)\|+|r-s|\\\notag
&\leq \left(1-\gamma_0+\eps(N+1)\right)\|p-q\|.
\end{align}
The assumed bound on $\eps$ implies $\eps<\frac{\gamma_0}{N+1}$ by Proposition 2(b), from which it follows $1-\gamma_0+\eps(N+1)<1$. This implies $\|H^j(p)-M_\eps\|\rw 0$ exponentially fast as $j\rw\infty$. \ $_\square$

\vspace{0.15in}
This completes the proof of Theorem 1. We now return to the model and use Theorem 1 to gain insight into the behavior of orbits of $H$ and, in particular, to show  the existence of a stable fixed point for which the ice line lies in tropical latitudes (the Jormungand climate state). We also   investigate system  bifurcations.

\section{Model behavior}

\subsection{Existence of the Jormungand climate state}
Consider the attracting invariant  manifold $M_\eps=\mbox{graph}(g^*)$ for mapping \eqref{H}, whose existence follows from Theorem 1. Using the fact $M_\eps$ is parametrized by $\eta$ and lies within $O(\eps)$ of $M_0$, the ice line equation on $M_\eps$ becomes
\begin{align}\notag
p_\eps(\eta)=\eta+G(g^*(\eta),\eta)&=\eta+\eps\left(\tx{\sum^N_{i=0}}(g^*(\eta))_i \, q_{2i}(\eta)-T_c\right)\\\notag
&=\eta+\eps\left(\tx{\sum^N_{i=0}}((h^0(\eta))_i+O(\eps)) \, q_{2i}(\eta)-T_c\right)\\\notag
&=\eta+\eps\left(\tx{\sum^N_{i=0}}f_{2i}(\eta) \, q_{2i}(\eta)-T_c\right)+O(\eps^2)
\end{align}
(noting $\sum^N_{i=0} q_{2i}(\eta)$ is bounded for $\eta\in\R$). Hence the dynamics of $H(x,\eta)$ on $M_\eps$ are well-approximated by  orbits of the function
\begin{equation}\label{p}
\phi_\eps:\R\rw\R, \ \phi_\eps(\eta)=\eta+\eps\left(\tx{\sum^N_{i=0}}f_{2i}(\eta) \, q_{2i}(\eta)-T_c\right),
\end{equation}
which we write
\begin{equation}\label{zee}
\phi_\eps(\eta)=\eta+\eps z(\eta), \ z(\eta)=\tx{\sum^N_{i=0}}f_{2i}(\eta) \, q_{2i}(\eta)-T_c.
\end{equation}
Iteration of \eqref{p} 
thus allows for the analysis of the behavior of system \eqref{extsystem} for $\eps$ positive and sufficiently small. 

We pause to comment on the choices for parameters presented in Table 1 and used for the plots below. Generally, parameters utilized  in a conceptual model of a system as complex as planetary climate will not be well-constrained. As mentioned previously, the parameters $A$ and $B$ used in the outgoing longwave radiation (OLR) term have been estimated as functions of surface temperature using satellite data \cite{graves}. During the extremely cold climate of the Neoproterozoic Era ice ages under consideration here, the drawdown of atmospheric CO$_2$ via silicate weathering would be greatly reduced. The associated buildup of
atmospheric CO$_2$ would serve to decrease OLR, leading to a smaller value for $A$ relative to that provided in \cite{graves}   for the present climate. The value $A=164$ Wm$^{-2}$ in Table 1 is in line with values used in \cite{abb}.

The value  $B=1.9$ W$(^\circ$C m$^2)^{-1}$ is taken from \cite{graves}. The critical temperature $T_c=0^\circ$C is from \cite{abb}, while $R=20$ Wyr($^\circ$C  m$^2$)$^{-1}$ is taken from \cite{schwartz}. The solar constant is known to have been 94\% of its current value 600-700   mya; hence we use $Q=321$ Wm$^{-2}$. 

From a conceptual modeling perspective, the diffusion constant appropriate for the present climate is $D\approx 0.5$ W/($^\circ$C  m$^2$) \cite{rayneo}. Given that the cold climates discussed here may reduce the efficiency of the heat transport, we use  $D\approx 0.25$, as indicated in \cite{rayneo}. Assuming the obliquity was not significantly different  600-700  mya, we use the current value $\beta=23.4^\circ$, for which the coefficients $s_{2i}, \ i=0, 1, ... , 5$ in Table 1  and arising from expansion \eqref{ess} are provided in \cite{dickalice}.

The albedo values are difficult to estimate (\cite{raybook}, chapter 3). The choice of the albedo of bare ice $\al_i=0.4$ being slightly  greater than that of ice free surface ($\al_1=0.3$) is important; for sufficiently large $\al_i$-values the ice albedo feedback dominates, resulting in the loss of the Jormungand solution described below. 

With parameters as in Table 1, we plot $z(\eta)$ for $N=1, ... , 5$ in Figure 3. The function $\phi_\eps(\eta)$ will have a fixed point $\eta=\eta^*$ precisely when $z(\eta^*)=0$. Moreover, any zeros $\eta=\eta^*$ of $z(\eta)$ satisfying $z^\pr(\eta^*)<0$ will correspond to attracting fixed points for $\phi_\eps(\eta)$, assuming $\eps>0$ is sufficiently small, since $\phi^\pr_\eps(\eta^*)=1+\eps z^\pr(\eta^*)$. Hence model equations \eqref{extsystem} admit an attracting fixed point having $\eta$-value near any $\eta=\eta^*$  for which  $z(\eta^*)=0$ and $z^\pr(\eta^*)<0$. Near any zero $\eta=\eta^*$ of $z(\eta)$ with $z^\pr(\eta^*)>0$, system \eqref{extsystem} will have an unstable fixed point with 1-dimensional unstable manifold  and $(N+1)$-dimensional stable manifold. 

\begin{figure}[t]
\includegraphics[width=5.5in,trim = 1.4in 5.53in 1.5in  1.2in, clip]{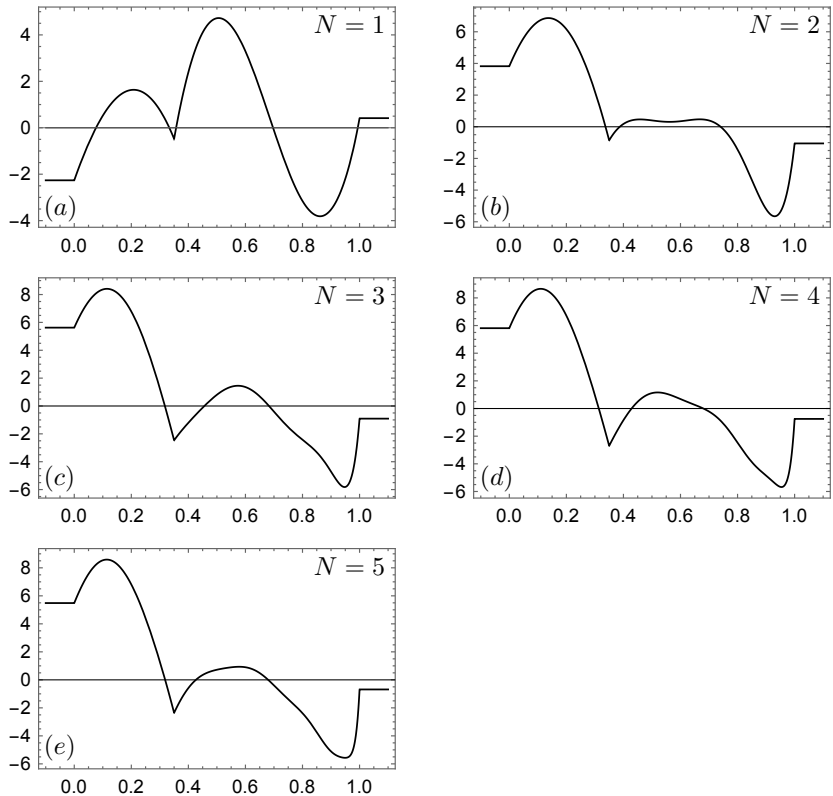}
\begin{center}
\caption{Plots of the function $z(\eta)$ \eqref{zee}. To each zero $\eta=\eta^*$  of $z(\eta)$  there corresponds a fixed point of $H(x,\eta)$ on the invariant manifold $M_\eps$. If $z^\pr(\eta^*)<0$ (resp., $z^\pr(\eta^*)>0$), the corresponding fixed point of $H$ on $M_\eps$ is stable (resp., unstable).
}
 \end{center}
\end{figure}

As for the dynamics on $M_\eps$ more globally, note $p_\eps$ is injective due to Proposition 3. Given that $g^*$ and each $q_{2i}$ are bounded for $\eta\in\R, \ p_\eps$ is surjective as well. Thus $p_\eps:\R\rw\R$ is a homeomorphism, from which it follows  all $p_\eps$-orbits are monotonic. For $\eta$ satisfying  $G(g^*(\eta),\eta)<0$, \  $p_\eps(\eta)<\eta$, and the corresponding orbit will be decreasing.  For $\eta$ satisfying  $G(g^*(\eta),\eta)>0$, \  $p_\eps(\eta)>\eta$, and the corresponding orbit will be increasing. 

For the parameter values in Table 1, $p_\eps$ will have three fixed points $\eta_1<\eta_2<\eta_3$ near the three respective zeroes of $z(\eta)$ for $N>1$ (see Figure 3(b-e)). Starting with an $\eta>\eta_2$, the ice line will tend to the stable fixed point with a smaller ice cap, for which $\eta=\eta_3$. For initial conditions with $\eta<\eta_2$, the ice line will tend to the stable fixed point with a large ice cap  ($\eta=\eta_1$). We note, in particular, that the model admits a stable Jormungand equilibrium state with the ice line resting in tropical latitudes below $y=\rho$ for each of $N=1, ... , 5$.
We also remark that both  an ice free Earth and a snowball Earth are unstable in the sense that $p_\eps(1)<0$ (so an ice line would descend to $\eta=\eta_3$), and   $p_\eps(0)>0$ (so the ice line would retreat to $\eta=\eta_1$), for $N>1$.

\subsection{Local bifurcations}

We fix $N=5$ in this section, for which $z(\eta)$ is plotted in Figure 3(e). We note the parameter $A$ can be thought of as a proxy for atmospheric greenhouse gas concentrations, in the following  sense. A large $A$-value corresponds a large OLR term, which in turn corresponds to a low value of, say, atmospheric CO$_2$ concentration. A smaller $A$-value leads to a small OLR term, and hence a larger 
atmospheric  concentration of CO$_2$.

Note the parameter $A$ appears in $f_{2i}(\eta)$ only when $i=0$. Hence $A$ appears precisely once in   expansion \eqref{zee} for $z(\eta)$. While we omit the details, setting $z(\eta)=0$ and solving for $A$ yields a correspondence between fixed points of $\phi_\eps$ and $A$, that is, between the position of the ice line at equilibrium and the parameter $A$.

For $\eta\in [0,1]$ the position of the ice line at equilibrium as a function of $A$ is plotted  in Figure 4. The solid branches in Figure 4 represent stable fixed points, while dotted branches are comprised of unstable fixed points. The vertical red line in Figure 4 corresponds to the $A$-value used in Figure 3, for which one unstable fixed point was sandwiched between two stable fixed points.

For sufficiently small $A$ values---or large CO$_2$ concentrations---there are no equilibria, with the ice line increasing for all time. The physical interpretation is that the climate is heading to an ice free state. A saddle-node bifurcation occurs at $A\approx 153$ Wm$^{-2}$, with a small unstable ice cap and a moderately sized stable ice cap equilibria pair appearing. The existence of the unstable small ice cap equilibrium state  was previously shown to occur when using Legendre polynomial expansions to approximate the temperature equation \eqref{Tdiff} in  \cite{north84}.

The stable Jormungand state comes into existence as $A$ passes through a nonsmooth fold bifurcation \cite{budd} at roughly 159 Wm$^{-2}$. Note the system now admits four fixed points, provided by two pair of stable/unstable fixed points. The stable moderately sized ice cap equilibrium disappears in  a saddle node bifurcation at $A\approx 166$ Wm$^{-2}$.
As the atmospheric CO$_2$ concentrations continue to decrease---and the OLR term increases, cooling the planet---the  Jormungand state is the sole remaining stable equilibrium. 

The model exhibits an unstable equilibrium with ice line near the equator for $A$ slightly larger than  175 Wm$^{-2}$. There is one final saddle-node bifurcation
at $A\approx 181$ Wm$^{-2}$ in which the Jormungand state is lost; the OLR is sufficiently large so as to dominate the warming effect the bare ice had on the absorption of insolation. If $A$ is large enough there are no equilibria, with the ice line heading toward the equator in a runaway snowball Earth episode.

\begin{figure}[t]
\includegraphics[width=6.1in,trim = 1.3in 7.35in .5in  1.25in, clip]{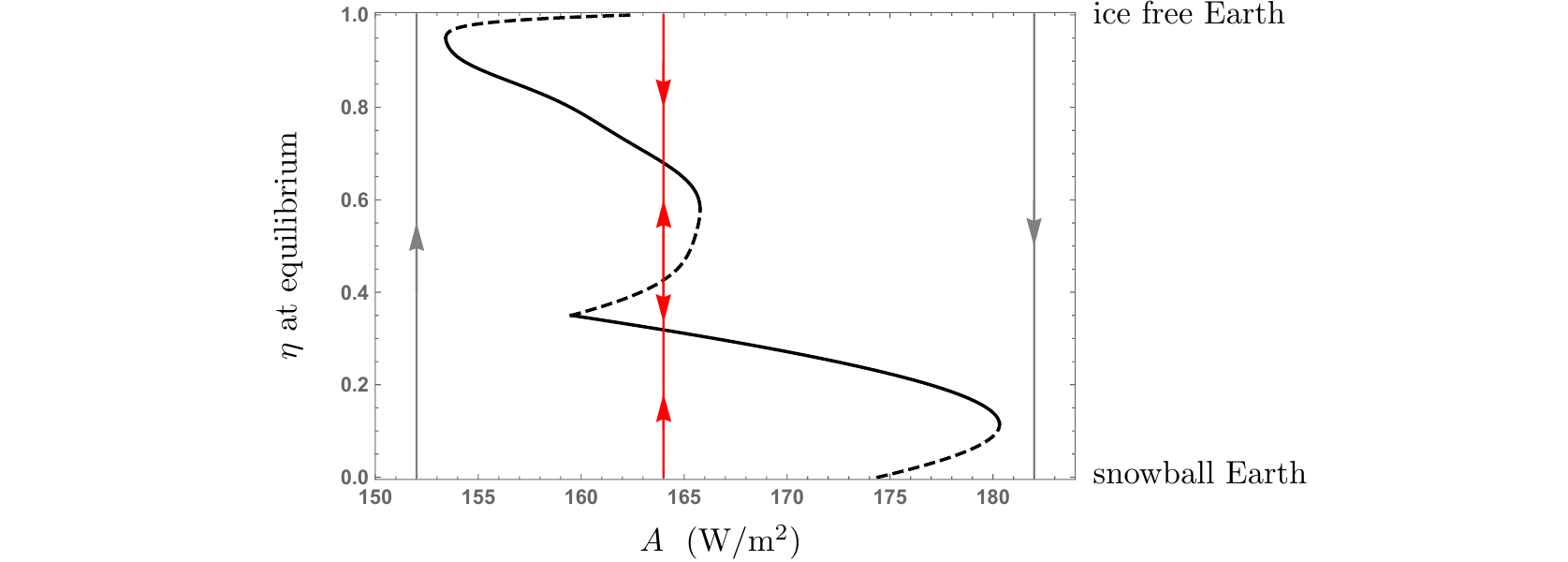}
\begin{center}
\caption{The bifurcation diagram as the parameter $A$ is varied, with $D=0.25$ W($^\circ$C m$^2)^{-1}$. Solid branches are stable fixed points and dashed branches are unstable fixed points. The red phase line represents the $A$-value used in Figure 3.
}
 \end{center}
\end{figure}

We note the range of $A$-values for which a stable Jormungand state exists is as large as that obtained in \cite{abb} when using   heat transport formulation \eqref{mean}. That the Jormungand state is stable for such a wide range of CO$_2$ values is notable in that the Neoproterozoic Era ice ages under consideration here lasted  several million years, evidently exhibiting an insensitivity to changes in CO$_2$.

Fixing $A=164$ Wm$^{-2}$, we plot the equilibrium ice line position as a function of the diffusion coefficient $D$ in Figure 5. While there is a significant range of $D$-values for which the stable Jormungand state exists, the model exhibits bistability over much of this range. The Jormungand state is the sole stable equilibrium only for $D$-values roughly between 0.35 and 0.44. For very efficient heat transport (to the right in Figure 5), the climate system tends toward a ``uniform" climate of either a warm ice free world or a snowball Earth, depending on whether the ice line starts above or below an unstable equilibrium position.
For highly inefficient heat transport (toward the left in Figure 5), the warm tropical climate and the cold northern regions have little interaction, leading to a stable equilibrium solution with the ice line lying between 0.4 and 0.6. 

\begin{figure}[t]
\includegraphics[width=6.1in,trim = 1.3in 7.35in .5in  1.25in, clip]{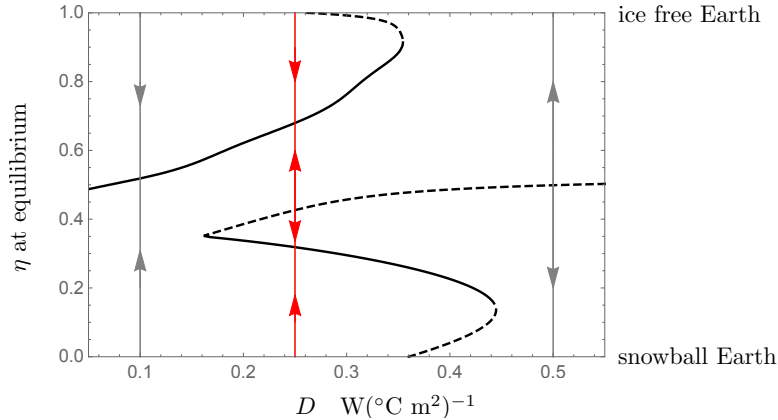}
\begin{center}
\caption{The bifurcation diagram as the parameter $D$ is varied, with $A=164$ Wm$^{-2}$. Solid branches are stable fixed points and dashed branches are unstable fixed points. The red phase line represents the $D$-value used in Figure 3.
}
 \end{center}
\end{figure}

\section{Conclusion}

M. Budyko's zonally averaged surface temperature model is an energy balance model comprised of  terms representing incoming solar radiation, outgoing longwave radiation, and meridional heat transport. Budyko introduced his conceptual model to investigate the role played by positive ice albedo feedback in influencing climate. E. Widiasih brought a dynamical systems approach to Budyko's model by coupling the temperature profile equation to an ice line allowed to respond to changes in temperature. Working with Budyko's heat transport term \eqref{mean} and in the infinite dimensional setting, with parameters aligned with the present climate, Widiasih showed there exists a stable equilibrium temperature function--ice line pair $(T^*(y),\eta^*)$, with $y=\eta^*$ positioned roughly as it is today. She also showed there exists a second, unstable equilibrium solution with a large ice cap. While not equilibria, Widiasih argued that a snowball Earth state is stable, while the ice free state is unstable.

The heat transport term \eqref{mean} used by Budyko was replaced  with the diffusive transport term \eqref{diffuse},
and the resulting equation coupled to Widiasih's ice line equation, in \cite{mediff}. Using the spectral method, and with parameters aligned with the present climate, results mirroring those of Widiasih were obtained.

Controversy continues to shroud the proposition that snowball Earth events have occurred in our planet's past. A means by which the effect of positive ice albedo feedback might be lessened for sufficiently large ice sheets, so that a strip of ocean water about the equator remained ice free,  was presented in \cite{abb}. Incorporating this mechanism into the Budyko-Widiasih model with heat transport term \eqref{mean}, and working in the infinite dimensional setting, a stable Jormungand equilibrium solution with the ice line in tropical latitudes was shown to exist in \cite{mewid}. It is of interest to note a stable Jormungand state was found when running an idealized general circulation model in \cite{abb}, parametrized for the extremely cold world associated with an extensive ice age.

We have incorporated the diffusive heat transport term \eqref{diffuse}, along with changes in surface albedo motivated by \cite{abb}, into the Budyko-Widiasih model in this work. Assuming a discrete time approach and using the spectral method, the analysis results in a nonsmooth singular perturbation problem. We proved an invariant manifold persists under perturbation, assuming the ice sheet moves sufficiently slowly relative to the evolution of temperature. This result was used to show the model exhibits a stable Jormungand equilibrium state, with appropriately chosen parameters. The model displays a rich bifurcation structure, both as atmospheric CO$_2$ concentrations vary and with changes in the efficiency of the meridional heat transport.

The dynamics of the model analyzed here are determined by a single equation governing the movement of the ice line, for $\eps>0$ and sufficiently small.  Coupling this one equation with either a dynamic equation for atmospheric CO$_2$ concentrations (as in \cite{anna}, where $A$ varies with time), or with a latitudinal dependence assumed for the diffusion constant (as suggested in \cite{north}), would  provide for   further lines of inquiry.

\vspace{0.3in}
\noindent
{\bf Acknowledgment.} The author recognizes and appreciates the support of the Mathematics and Climate Research Network (www.mathclimate.org).

\renewcommand{\refname}

\end{document}